\documentclass[10pt,twoside]{article}
 \usepackage{amssymb,latexsym,amsmath,amsthm}
 \usepackage[numbers,sort&compress]{natbib}

 \usepackage{fullpage}

\newtheorem{thm}{Theorem}[section]
\newtheorem{lem}{Lemma}[section]
\newtheorem{prop}{Proposition}[section]

\newtheorem{remark}{Remark}
\newtheorem{dfn}{Definition}
\newtheorem{notation}{Notation}

\begin{document}
\date{}

\title{Effect of weights on stable solutions of a quasilinear elliptic equation}
\maketitle

\begin{center}
\author{{Mostafa Fazly\footnote{Current affiliation:  Department of Mathematical and Statistical Sciences, University of Alberta, Edmonton, Alberta, T6G 2G1, Canada}}
}
\\
{\it\small Department of Mathematics}\\ {\it\small University of British Columbia}, 
{\it\small Vancouver, B.C. Canada V6T 1Z2}\\
{\it\small e-mail: fazly@math.ubc.ca}\vspace{1mm}\\

\end{center}

\vspace{3mm}

\begin{abstract}
In this note, we study Liouville theorems for the stable and finite Morse index weak solutions of the quasilinear elliptic equation $-\Delta_p u= f(x) F(u) $ in $\mathbb{R}^n$ where $p\ge 2$, $0\le f\in C(\mathbb{R}^n)$ and $F\in C^1(\mathbb{R})$.  We refer to $f(x)$ as {\it weight} and to  $F(u)$ as {\it nonlinearity}.  The remarkable fact is that if the weight function is bounded from below by a strict positive constant that is $0<C\le f$ then it does not have much impact on the stable solutions, however, a nonnegative weight that is $0\le f$ will push certain critical dimensions.  This analytical observation has potential to be applied in various models to push certain well-known critical dimensions. 

For a general nonlinearity $F\in C^1(\mathbb{R})$ and $f(x)=|x|^\alpha$, we prove Liouville theorems in dimensions $n\le \frac{4(p+\alpha)}{p-1}+p$, for bounded radial stable solutions. For specific nonlinearities $F(u)=e^u$, $u^q$ where $q>p-1$ and $-u^{q}$ where $q<0$,  known as the Gelfand, the Lane-Emden and the negative exponent nonlinearities, respectively, we prove Liouville theorems for both radial finite Morse index (not necessarily bounded) and stable (not necessarily radial nor bounded) solutions.  

\end{abstract}


\section{Introduction} 
We examine weak solutions of the following equation 
\begin{equation} \label{gener} 
-\Delta_p u= f(x) F(u) \ \  \ \text{in}\ \ \mathbb{R}^n,
\end{equation}
 where $p\ge2$, $F\in C^1(\mathbb{R})$ and $f\in C(\mathbb{R}^n)$.   The notation $\Delta_p$ stands for the $p$-Laplacian operator given by $\Delta_p u=\text{div} (|\nabla u|^{p-2}\nabla u)$.  This is a generalization of the following equation
 \begin{equation}\label{p=2gener} 
-\Delta u= F(u) \ \  \ \text{in}\ \ \mathbb{R}^n,
\end{equation}
that has been of interest for various nonlinearities $F$. For sign changing nonlinearities this equation is of great interest for the double well potential nonlinearities, e.g., $F(u)=u-u^3$ known as the Allen-Cahn equation.  For nonsign changing  nonlinearities this equation is of interest for various nonlinearities such as  $F(u)=e^u$, $F(u)=u^q$ where $q>p-1$ and $F(u)=-u^{q}$ where $q<0$,  known as the Gelfand, the Lane-Emden and the negative exponent nonlinearities, respectively. Note that  when the nonlinearity is the negative exponent nonlinearity $F(u)=-u^{-2}$, (\ref{p=2gener}) is the counterpart of the second order elliptic MEMS equation on the whole space $ \mathbb{R}^n$.   The second order elliptic MEMS equation on a bounded domain $\Omega\subset \mathbb{R}^n$ is the following equation
\begin{equation*}
\left\{
                      \begin{array}{ll}
                        -\Delta u=(1-u)^{q}  & \hbox{in $\Omega$;} \\
                       u=0 & \hbox{on $\partial\Omega$.} \\
                           \end{array}
                    \right.
                    \end{equation*}  
where  $q<0$. The MEMS equation arises in modelling an electrostatic Mircro-ElectroMechanical System (MEMS) device. For details we refer the interested readers to the book of Pelesko-Bernstein \cite{pb} for physical derivations of the MEMS model and to the book of Esposito-Ghoussoub-Guo \cite{egg2} for the mathematical analysis of the model. See also the interesting papers \cite{gpw,gg}. One of the primary goals in the design of MEMS devices is to
optimize the pull-in distance over a certain allowable voltage range that is set by the
power supply. Here the pull-in distance refers to the maximum stable deflection of
the elastic membrane before quenching occurs.

The weight $f(x)\ge 0$ makes the equation (\ref{gener}) much more challenging and as a general statement, some standard techniques such as moving plane methods and Sobolev embeddings cannot be applied anymore.  In particular, the power weight $f(x)=|x|^\alpha$ has been of interest in this context and it was introduced by M. H\'{e}non \cite{h} in equation 
\begin{equation}\label{henon}
\left\{
                      \begin{array}{ll}
                       -\Delta u=|x|^\alpha u^{q}  & \hbox{for  $|x|<1$,} \\
                       u=0 & \hbox{on $|x|=1$} \\
                           \end{array}
                    \right.
                    \end{equation}  
to model and study spherically symmetric clusters of stars.  This equation is now known as the H\'{e}non equation for $\alpha>0$ and the H\'{e}non-Hardy equation for $\alpha<0$.  Ten years later, Ni in \cite{n} explored properties of positive radial solutions of the H\'{e}non equation on the unit ball and observed the fact that the power profile $f(x)=|x|^\alpha$ enlarges considerably the range of solvability beyond the classical critical threshold, i.e., $q<2^*-1=\frac{n+2}{n-2}$ to $q<2^*-1+\frac{2\alpha}{n-2}=\frac{n+2+2\alpha}{n-2}$ where $2^*$ is the critical exponent for the Sobolev embedding $H_0^1\hookrightarrow L^q$.   On the other hand, as it is shown by Smets, Su and  Willem in \cite{ssw} and references therein, equation (\ref{henon}) also admits nonradial solutions for $q<2^*-1$.  The existence of nonradial solutions for the full range $q<\frac{n+2+2\alpha}{n-2}$ is still an open problem.     Note that since the function $|x|\to |x|^\alpha$ is increasing, the classical moving planes arguments given by Gidas, Ni and Nirenberg   in \cite{gnn} cannot be applied to prove the radial symmetry of the solutions of (\ref{henon}). Therefore, the existence of nonradial solutions for this equation is natural.

The effect of power-like permittivity profiles $f(x)=|x|^\alpha$ and the dimension $n$ on both the pull-in voltage and the pull-in distance of the MEMS equation has been studied in \cite{cg,gpw}. Also the effect of the power weights on the semilinear elliptic equations with various nonlinearities and in the notion of stability has been  extensively studied for both bounded and unbounded domains in \cite{e1,e2,egg1,egg2,gpw,cg,wy,cf,mostafagh,mostafa}.  Our motivation to write this note is to study the impact of the weight functions and specially power weights on stable solutions of the weighted quasilinear equation (\ref{gener}) for various nonlinearities.     Quasilinear equations and models play a fundamental role in mathematical analysis and applied sciences, e.g. the minimal surface equation, various models of macroscopic conservation laws in gas dynamics and quasilinear reaction-diffusion equations.

 In what follows, we show that strictly positive weights do not change the critical dimensions for the nonexistence of stable solutions, see Remark \ref{remarkM}. However, nonnegative weights push certain critical dimensions. More precisely, we first explore properties of radial solutions of (\ref{gener}) with a general nonlinearity $F\in C^1(\mathbb{R})$. We  prove Liouville theorems for bounded stable solutions, following ideas given in \cite{cc,ces,sv} and using a suitable change of variable to eliminate weights.  Then, applying the Harnack's inequality, we prove Liouville theorems for finite Morse index solutions of (\ref{gener}) with three different nonlinearities as well as some Liouville theorems for not necessarily radial nor bounded stable solutions. For nonradial solutions, we apply the Moser iterations method that is multiplying  equation (\ref{gener}) by powers of the nonlinearity $F(u)$ as well as certain integral estimates.  The Moser iterations method in this context was developed in \cite{egg1,egg2,f1,f2}.

Throughout the paper,  we assume that $\alpha\ge0$, $n>p \ge 2$ and solutions of (\ref{gener}) are positive for power nonlinearities. However, most of proofs can be adapted for $\alpha+p>0$ and for sign changing solutions of (\ref{gener}) if $F(u)=|u|^{q-1}u$ for $q>p-1$. Since $\Delta_p$ has variational structure, we can associate the following energy functional to (\ref{gener}),
\begin{equation*}
I(u):=\int_{\mathbb{R}^n} \frac{1}{p} |\nabla u|^{p}-f(x) \mathcal{F}(u),
\end{equation*}
 where $\mathcal{F}(t)=\int_{0}^{t}F(s) ds$.  Here we have some definitions.
\begin{dfn}\label{def}
We call $u\in C^{1,\gamma} {(\mathbb{R}^n)}$
\begin{itemize}
\item  {\it weak solution} of (\ref{gener}), if for all $\phi\in C_c^1(\mathbb{R}^n)$ the following holds: 
 \begin{equation} \label{wgener} 
\int_{\mathbb{R}^n} |\nabla u|^{p-2} \nabla u\cdot \nabla \phi = \int_{\mathbb{R}^n} f(x)F(u) \phi.
\end{equation} 

  \item {\it stable} solution of (\ref{gener}) if not only $u$ is a solution of (\ref{gener}), but also the second variation of the associated energy functional is nonnegative. This means for $\phi\in C_c^1(\mathbb{R}^n)$, we have
 $$I_{uu}(\phi):=\int_{\mathbb{R}^n} |\nabla u|^{p-2} |\nabla \phi|^2 + (p-2) \int_{\mathbb{R}^n} |\nabla u|^{p-4} |\nabla u\cdot\nabla \phi|^2 - \int_{\mathbb{R}^n} f(x)F'(u)\phi^2\ge 0. $$
In particular, for $p\ge 2$ we have 
\begin{equation} \label{stableineq} 
\int_{\mathbb{R}^n}  f(x)  F'(u) \phi^2 \le (p-1) \int_{\mathbb{R}^n} | \nabla u|^{p-2}|\nabla\phi|^2 \qquad \forall \phi \in C_c^1(\mathbb{R}^n).
\end{equation}
\item  {\it stable outside a compact set} $\Sigma\subset\mathbb{R}^n$ if $I_{uu}(\phi)\ge 0$ for all $\phi\in C_c^1(\mathbb{R}^n\setminus\Sigma)$.  Also we say $u$ has a {\it Morse index} equal to $m\ge 1$ if $m$ is the maximal dimension of a subspace $X_m$ of $C^1_c(\mathbb{R}^n)$ such that $I_{uu}(\phi)<0$ for all $\phi\in X_m\setminus \{0\}$.
\end{itemize}
\end{dfn}

Note that if $u$ is of Morse index $m$, then there exist $\phi_1, ..., \phi_m$ such that $X_m= span\{\phi_1, ..., \phi_m\}\subset C_c^1(\mathbb{R}^n)$ and $I_{uu}(\phi)<0$ for all $\phi\in X_m\setminus \{0\}$. So, for all $\phi\in C_c^1(\mathbb{R}^n\setminus \Sigma)$ we have $I_{uu}(\phi)\ge 0$, where $\Sigma = \cup_{i=1}^{m} supp(\phi_i)$. Therefore, $u$ is stable outside a compact set $\Sigma\subset \mathbb{R}^n$. On the other hand, the $p$-Laplacian operator of a radial function in dimension $n$ is given  by $\Delta_{p,n} u=|u_r|^{p-2}\left( (p-1)u_{rr} +\frac{n-1}{r} u_r\right).$
 The operator  $\Delta_{p,n} $ depends on $p$ and $n$ and it is well-defined for any $p,n\in\mathbb{R}$. From Definition \ref{def}, a radial weak solution $u$ of (\ref{gener}) satisfies 
 \begin{equation} \label{wrgener} 
\int_{0}^{\infty} r^{n-1}| u_r|^{p-2} u_r \phi_r = \int_{0}^{\infty} r^{n-1} f(r)F(u) \phi, \ \ \ \ \ \ \text {for all} \ \  \phi\in C_c^1\left([0,\infty)\right).
\end{equation}

\begin{notation} For the sake of simplicity in computations and presentations, we define the following notations. 
\begin{itemize} 
\item $N_p(\alpha):=\frac{1}{p}  \left( 1+\frac{\alpha}{p}  \right)  \left(p+2-\frac{p(N+\alpha)}{p+\alpha}+2\sqrt{   \frac{1}{p+\alpha} \left(   \frac{p(N-1)}{p-1}+\alpha\right)}\right)$.

\item  $q_p^*(\alpha):=\frac{p(q+\alpha+1)}{q-p+1} $ \ and \  $q_p^\#(\alpha):=\frac{p(q-1)+\alpha(p-2)}{q-p+1}$. 

\item   $q^+_p(\alpha):= \frac{q-1}{q-p+1}p+\frac{p-2}{q-p+1}\alpha+2(p+\alpha)\frac{q+\sqrt{q(q-p+1)}}{(p-1)(q-p+1)}$.
\item  $q^-_p(\alpha):=\frac{q-1}{q-p+1}p+\frac{p-2}{q-p+1}\alpha+2(p+\alpha)\frac{q-\sqrt{q(q-p+1)}}{(p-1)(q-p+1)}$.
\end{itemize}
\end{notation}

\section{Main Results}

The interesting fact about radial solutions of (\ref{gener}) is that by an appropriate change of variable, the  equation (\ref{gener}) with the power weight $f(x)=|x|^\alpha$ in dimension $n$ can be modified as an equation with a constant weight and in a new fractional dimension. To get this change of variable, we have used ideas given by  Cowan-Ghoussoub in \cite{cg} for the semilinear equations, i.e. $p=2$. They applied  this change of variable on the unit ball to provide proofs for various phenomena observed by Guo-Pan-Ward in \cite{gpw}. This change of variable is also used in \cite{bg} for the semilinear equations and on the whole space.

\begin {prop} \label{change} 
The radial  function $u(r)$ satisfies $-\Delta_{p,n} u(r)= 
r^\alpha F(u(r))$ in $B_{R}$ in dimension $n$ if and only if the function $\omega(s):=u(r) $ where $s:=r^{1+\frac{\alpha}{p}}$ is a solution of $-\Delta_{p,\frac{p(\alpha+n)}{\alpha+p}} \omega(s)= (1+\frac{\alpha}{p})^{-p} F(\omega(s))$ in  $B_{R^{  1+\frac{\alpha}{p} }}$ in the fractional dimension $\frac{p(n+\alpha)}{p+\alpha}$. Also, a similar result holds on the whole space $\mathbb{R}^{n}$.

\end{prop}

\begin{remark}
By scaling one can remove the constant $(1+\frac{\alpha}{p})^{-p}$. 
\end{remark}

\noindent \textbf{Proof:}  Set $\omega(s):=u(r)$ for $s:=r^{1+\frac{\alpha}{p}}$.  Then, by a straightforward calculation, we have 
\begin{eqnarray*}
|u_r|^{p-2} \left(    (p-1) u_{rr} +\frac{n-1}{r} u_r     \right) =\left(1+\frac{\alpha}{p}\right)^p r^\alpha |\omega_s|^{p-2}  \left(    (p-1) \omega_{ss} +\frac{N_{\alpha,p}-1}{s} \omega_s     \right) 
\end{eqnarray*}
where $N_{\alpha,p}:=\frac{p(n+\alpha)}{p+\alpha}$. Since the $p$-Laplacian operator is given by $\Delta_{p,n} u= |u_r|^{p-2} \left(    (p-1) u_{rr} +\frac{n-1}{r} u_r     \right)$, we get 
$$\Delta_{p,n} u(r)=\left(1+\frac{\alpha}{p}\right)^p r^\alpha \Delta_{p,\frac{p(n+\alpha)}{p+\alpha}} \omega(s).
$$
Therefore,    
$$-\Delta_{p,\frac{p(n+\alpha)}{p+\alpha}} \omega=\left(1+\frac{\alpha}{p}\right)^{-p} F(\omega).
$$
\hfill $\Box$

In the following theorem, we prove a Liouville theorem for radial solutions of quasilinear equation (\ref{gener}) with a general nonlinearity $F\in C^1(\mathbb{R})$. Our methods of proofs are the methods developed by Cabr\'e and Capella \cite{cc} for the case $p=2$ and $\alpha=0$.

\begin{thm} \label{radialgener2}
Let $n> 1-\alpha+\frac{\alpha}{p}$, $F\in C^1(\mathbb{R})$ and $u$ be a nonconstant bounded radial stable weak solution of (\ref{gener}). Then,   $ n > \frac{4(p+\alpha)}{p-1}+p$ and $$|u(r)-u_{\infty}|\ge C r^{N_p(\alpha)},  \ \ \ \forall r\ge 1, $$
where $u_{\infty}=\lim_{r\to\infty} u(r)$ and $C$ does not depend on $r$. 
\end{thm}

Moreover, we have the following pointwise estimate for not necessarily bounded solutions of (\ref{gener}).

\begin{thm} \label{radialgenerest}
Let $n> 1-\alpha+\frac{\alpha}{p}$, $F\in C^1(\mathbb{R})$ and $u$ be a nonconstant radial stable weak solution of (\ref{gener}). Then,   there exist positive constants $C$ and $r_0$ such that for $r\ge r_0$ we have 
\begin{equation}\label{radialestimate}
|u(r)|\ge  C \left\{
                      \begin{array}{ll}
                        r^{N_p(\alpha)}, & \hbox{if $n\neq \frac{4(p+\alpha)}{p-1}+p$;} \\
                       \ln (r), & \hbox{if $n=\frac{4(p+\alpha)}{p-1}+p$.} \\
                           \end{array}
                    \right.
                    \end{equation}     
  The constant $C$ does not depend on $r$.                 
\end{thm}

Applying the Harnack's inequality, we show the following Liouville theorem for radial solutions of (\ref{gener}) of finite Morse index. Note that there is no boundedness assumption in this theorem.

\begin{thm}  \label{radialcompact}
There is no radial weak solution with finite Morse index  of (\ref{gener}) if
\begin{itemize}  
  \item[(i)]  $F(u)= e^u$ in dimensions $p<n<\frac{4(p+\alpha)}{p-1}+p$.

  \item[(ii)]    $F(u)= u^q$  with  $q>p-1$  in dimensions $q_p^*(\alpha)<n<q^+_p(\alpha)$.
  
  \item[(iii)]    $F(u)= -u^q$   with  $q<0$  in dimensions $q_p^\#(\alpha)<n<q^-_p(\alpha)$.
\end{itemize} 
\end{thm}

For $p=2$ and $\alpha=0$, the above results $(i), (ii), (iii)$ and in corresponding dimensions  $$2<n<10, \ \ \ \  2=q_2^*(0)<n<q^+_2(0) \ \ \ \text{and}\ \ \ \  2=q_2^\#(0)<n<q^-_2(0)$$ are given in \cite{df}, \cite{f1} and \cite{egg2}, respectively.  Note that the term $q_p^\#(\alpha)=\frac{p(q-1)+\alpha(p-2)}{q-p+1}$ that appears as the lower bound for the dimension in $(iii)$ is always less than or equal to $p$ if and only if $(p+\alpha)(p-2)\ge 0$. So, the lower bound in $(iii)$ covers the lower dimensions $p<n$.
On the other hand, the term $q_p^*(\alpha)=\frac{p(q+\alpha+1)}{q-p+1}$ that appears in $(ii)$ is a special exponent in a sense that applying the Pohozaev identity, Proposition \ref{change}  and the standard techniques given by Caffarelli-Gidas-Spruck in \cite{cgs} for radial solutions, one can see that there is no nontrivial radial weak solution of (\ref{gener}) if $F(u)=u^q$ for $q>p-1$ in dimensions 
$ n<q_p^*(\alpha).$

Equation (\ref{gener}) for $F(u)=u^q$ when $q>p-1=1$ and $f(x)=|x|^\alpha$ is a generalization of the Lane-Emden equation
\begin{equation}\label{lane}
-\Delta u=u^q\ \  \ \text{in}\ \ \mathbb{R}^n.
\end{equation}
The positive classical solutions of (\ref{lane}) have been completely classified. Gidas and Spruck in \cite{gs1,gs2} proved Liouville theorems on the whole space $\mathbb{R}^{n}$ and in the absence of stability for the Sobolev subcritical exponent $q<\frac{n+2}{n-2}$ or  equivalently $n<q_2^*(0)=\frac{2(q+1)}{q-1}$. They conjectured that there is no positive solution for the weighted equation
\begin{equation}\label{lanealpha}
-\Delta u=|x|^\alpha u^q\ \  \ \text{in}\ \ \mathbb{R}^n,
\end{equation}
 for  all $n<q_2^*(\alpha)$ and $\alpha>-2$.  Note that (\ref{lanealpha}) is a special case of  (\ref{gener}) for $F(u)=u^q$ when $q>p-1=1$ and $f(x)=|x|^\alpha$.  This conjecture is solved by Phan and Souplet in \cite{phs} when either $-2<\alpha<0$ and all dimensions or  $\alpha>0$ and in dimension $n=3$  for bounded solutions.  This problem has been open in higher dimensions.  We would like to mention that the author and Ghoussoub in \cite{mostafagh} have proved  a Liouville theorem for the finite Morse index solutions of (\ref{lanealpha}) for all dimensions $n\ge 3$ ae well as a Liouville theorem for bounded solutions of the corresponding system.

For the critical case $n=\frac{2(q+1)}{q-1}$,  Caffarelli, Gidas and Spruck \cite{cgs} proved that all solutions of (\ref{lane}) are given by $$ \frac{\lambda^{\frac{n-2}{2}}}{ ( 1+  \lambda^2 |x-x_0|^2  )^{\frac{n-2}{2}} }$$
where $\lambda>0$ and $x_0\in\mathbb{R}^n$ and this is extremal function for the well-known Sobolev inequality 
$$\int_{\mathbb{R}^n} |\nabla u|^2\ge S \left( \int_{\mathbb{R}^n} |u|^{\frac{2n}{n-2}  } \right)^{\frac{n-2}{n}}$$
Following the ideas of the mentioned results, we have the following Liouville theorem for finite Morse index solutions of  (\ref{gener})  for $F(u)=u^q$ and $q>p-1$ when dimension $n$ is either subcritical or critical.

\begin{thm} \label{poho1}
Let $u$ be a nonnegative weak solution with finite Morse index, not necessarily radial, of (\ref{gener})  for $F(u)=u^q$ when $q>p-1$. Then, $u$ is the trivial solution  in dimensions $ n<q_p^*(\alpha)$. In the critical dimension $n=q_p^*(\alpha)$ all radial solutions of (\ref{gener}) are of the following form 
\begin{equation}\label{radial}
u_{\epsilon}(r):= \left(\epsilon(n+\alpha)(\frac{n-p}{p-1})^{p-1}\right)^{\frac{n-p}{p(p+\alpha)}}(\epsilon+r^{\frac{p+\alpha}{p-1}})^{\frac{p-n}{p+\alpha}},
\end{equation}
where $\epsilon>0$ and they are stable outside a compact set $\overline{B_{R_0}}$ for an appropriate $R_0$. 
\end{thm}

For not necessary radial nor bounded solutions, a similar Liouville theorem can be proved as following.  Note that higher dimensions are the same as given in Theorem \ref{radialcompact}.

\begin{thm}  \label{stablenon}
There is no entire stable weak solution for (\ref{gener}) with one of the following nonlinearities 
\begin{itemize}  
  \item[(i)] $F(u)= e^u$ in dimensions $n<\frac{4(p+\alpha)}{p-1}+p$.

  \item[(ii)]   $F(u)= u^q$  with  $q>p-1$  in dimensions $n<q^+_p(\alpha)$. 
  
  \item[(iii)]  $F(u)= -u^q$   with  $q<0$  in dimensions $n<q^-_p(\alpha)$.
\end{itemize} 
\end{thm}

In short, the weight functions and in particular $f(x)=|x|^\alpha$ can be applied in certain models to push some known critical dimensions. This phenomenon is of interest in terms of mathematical analysis of models and also their applications.

%%%%%%%%%%%%%%%%%%%%%%%%%%%%%%%%%%%%%%%%%%%%%%
%%%%%%%%%%%%%%%%%%%%%%%%%%%%%%%%%%%%%

\section{Proofs and Ideas}

\subsection{Proofs of Theorem \ref{radialgener2} and \ref{radialgenerest}}\label{1}

The first Liouville theorem in the notion of stability for radial stable solutions of (\ref{gener}) with a general nonlinearity $F\in C^1(\mathbb{R})$, $f=1$ and $p=2$ was nicely proved by Cabr\'e and Capella \cite{cc}.  Then,  Castorina-Esposito-Sciunzi in \cite{ces} extended Cabr\'e-Capella results to the quasilinear equation (\ref{gener}) with $f=1$. Recently, Villegas \cite{sv} by deriving some technical estimates and using a new test function in the stability condition improved Cabr\'e-Capella's results for the semilinear case $p=2$ and $f=1$.  
\\
 To prove a Liouville theorem result for (\ref{gener}) with $f(x)=|x|^\alpha$, one can apply Cabr\'e-Capella ideas in \cite{cc}, Castorina-Esposito-Sciunzi results in \cite{ces} and Proposition \ref{change} to show that if $u$ is a nonconstant bounded radial weak solution of (\ref{gener}) in dimensions $ 1-\alpha+\frac{\alpha}{p} < n\le \frac{p+\alpha}{p(p-1)}(3p-1+2\sqrt{2p-1})-\alpha$, then $u$ is unstable. Also, the same result holds for $$ \frac{p+\alpha}{p(p-1)}(3p-1+2\sqrt{2p-1})-\alpha< n\le \frac{4(p+\alpha)}{p-1}+p,$$ if we assume that $\lim_{s\to s_0} |F'(s)||s-s_0|^{-q}=a\in (0,\infty)$, for every zero point $s_0$ of $F$ and for some $q=q(s_0)\ge0$. Note that this nondegeneracy condition on $F$ is satisfied if $F$ is a nonzero analytic function.  In what follows, we manage to drop the nondegeneracy assumption. To do so, we adapt some estimates given by Villegas in \cite{sv} for the nonweighted semilinear equation.  Here is one of our main integral estimates.
 
\begin{lem} \label{radialgenerlem1}
Let $n> 1-\alpha+\frac{\alpha}{p}$, $F\in C^1(\mathbb{R})$ and $u$ be a nonconstant radial stable weak solution of (\ref{gener}). Then,  for a given $S$ there exists $C>0$ such that 
\begin{equation}\label{estimate}
\int_{s}^{S} \frac{dt}{  t ^{\frac{p(n-1)+\alpha(p-1)}{p+\alpha}} {|\omega_s(t)|}^p  } \le C \ s^{-2\sqrt{   \frac{1}{p+\alpha} \left(   \frac{p(n-1)}{p-1}+\alpha\right)}},   \ \ \forall \ 1\le s\le S,
\end{equation}
where $\omega(s):=u(r)$ for $s=r^{1+\frac{\alpha}{p}}$. The constant $C$ does not depend on $s$ and $S$.
\end{lem}

\noindent\textbf{Proof:} Proposition \ref{change} guarantees $\omega(s)$ is a radial stable solution of 
\begin{eqnarray}
-\Delta_{p,\frac{p(n+\alpha)}{p+\alpha}} \omega (s)=\left(1+\frac{\alpha}{p}\right)^{-p} F(\omega (s))=:\tilde{F}(\omega(s)).
\end{eqnarray}
By the same argument as in the proof of Theorem 1 in \cite{cc}, Theorem 1.4 in \cite{ces} and Lemma 2.3 in \cite{sv}, one has
\begin{eqnarray} \label{stabilityradial}
(N_{\alpha,p}-1) \int_{0}^{\infty} t^{N_{\alpha,p}-3} |\omega_s(t)|^p \eta^2(t) dt \le (p-1) \int_{0}^{\infty} t^{N_{\alpha,p}-1} |\omega_s(t)|^p \eta^2_t(t) dt,
\end{eqnarray}
for all $\eta\in C_c^1[0,\infty)$, where $N_{\alpha,p}:=\frac{p(n+\alpha)}{p+\alpha}$. This is true for either bounded or unbounded solutions.

Now, set the following test function $\eta\in H^1(\mathbb{R}^n)\cap L^{\infty} (\mathbb{R}^n)$:
$$\eta(t):=   \left\{
                      \begin{array}{ll}
                        1, & \hbox{if $0\le t \le 1$;} \\
                       t^{-\sqrt{\frac{N_{\alpha,p}-1}{p-1}}}, & \hbox{if $1\le t \le s$;} \\
                        \frac{    s^{-\sqrt{\frac{N_{\alpha,p}-1}{p-1}}}   }{  \int_{s}^{S} \frac{dz}{z^{N_{\alpha,p}-1}|\omega_s(z)|^p  } }    \int_{t}^{S} \frac{dz}{z^{N_{\alpha,p}-1}  |\omega_s(z)|^p}  , & \hbox{if $s\le t \le S$;} \\
                          0, & \hbox{if $S\le t $.} \\
                           \end{array}
                    \right.$$           
By a straightforward calculation, for the given test function $\eta$, L.H.S. of (\ref{stabilityradial}) has the following lower bound,
\begin{eqnarray*}
(N_{\alpha,p}-1) \int_{0}^{\infty} t^{N_{\alpha,p}-3} |\omega_s(t)|^p \eta^2(t) dt &\ge& (N_{\alpha,p}-1) \int_{0}^{1} |\omega_s(t)|^p  t^{N_{\alpha,p}-3} dt \\
& &+ (N_{\alpha,p}-1) \int_{1}^{s} |\omega_s(t)|^p  t^{   -2 \sqrt{\frac{N_{\alpha,p}-1}{p-1}}     +N_{\alpha,p}-3} dt.
\end{eqnarray*}
On the other hand, since 
$$\eta_t(t)=   \left\{
                      \begin{array}{ll}
                        0, & \hbox{if $0\le t <1$;} \\
                       -\sqrt{\frac{N_{\alpha,p}-1}{p-1}} t^{-\sqrt{\frac{N_{\alpha,p}-1}{p-1}}-1}, & \hbox{if $1< t < s$;} \\
                       - \frac{    s^{-\sqrt{\frac{N_{\alpha,p}-1}{p-1}}}   }{  \int_{s}^{S} \frac{dz}{z^{N_{\alpha,p}-1}|\omega_s(z)|^p  } }    \frac{1}{t^{N_{\alpha,p}-1}  |\omega_s(t)|^p}  , & \hbox{if $s < t < S$;} \\
                          0, & \hbox{if $S< t $;} \\
                           \end{array}
                    \right.$$    
R.H.S. of (\ref{stabilityradial}), for the given $\eta$,  is the same as
\begin{eqnarray*}
(p-1) \int_{0}^{\infty} t^{N_{\alpha,p}-1} |\omega_s(t)|^p \eta^2_t(t) dt &=& (p-1)  \int_{1}^{s} \left(    \frac{  N_{\alpha,p}-1 }{p-1}\right)  t^{  -2\sqrt{\frac{N_{\alpha,p}-1}{p-1}} + N_{\alpha,p}-3} |\omega_s(t)|^p dt \\
& & +(p-1) \int_{s}^{S} t^{N_{\alpha,p}-1} |\omega_s(t)|^p     \frac{    s^{-2\sqrt{\frac{N_{\alpha,p}-1}{p-1}}}   } {  \left(\int_{s}^{S} \frac{dz}{z^{N_{\alpha,p}-1}|\omega_s(z)|^p  }\right)^2  }   \frac{dt} {\left(t^{N_{\alpha,p}-1}  |\omega_s(t)|^p\right)^2} \\
&= &(N_{\alpha,p}-1 ) \int_{1}^{s}   t^{  -2\sqrt{\frac{N_{\alpha,p}-1}{p-1}} + N_{\alpha,p}-3} |\omega_s(t)|^p dt \\
&&+ (p-1)  \frac{    s^{-2\sqrt{\frac{N_{\alpha,p}-1}{p-1}}}   } {  \int_{s}^{S} \frac{dz}{z^{N_{\alpha,p}-1}|\omega_s(z)|^p  }  } .
\end{eqnarray*}
Therefore, from (\ref{stabilityradial}) we obtain
\begin{eqnarray*}
(N_{\alpha,p}-1) \int_{0}^{1} t^{N_{\alpha,p}-3} |\omega_s(t)|^p dt\le (p-1)  \frac{    s^{-2\sqrt{\frac{N_{\alpha,p}-1}{p-1}}}   } {  \int_{s}^{S} \frac{dz}{z^{N_{\alpha,p}-1}|\omega_s(z)|^p  }  }  ,  \ \ \ \ \forall 1\le s\le S.
\end{eqnarray*}
 Hence, we have 
\begin{eqnarray*}
\int_{s}^{S} \frac{dz}{z^{N_{\alpha,p}-1}|\omega_s(z)|^p  } \le C \ s^{-2\sqrt{\frac{N_{\alpha,p}-1} {p-1}}} ,   \ \ \ \ \forall 1\le s\le S,
\end{eqnarray*}
where  $C:=\frac{p-1}{   (N_{\alpha,p}-1) \int_{0}^{1} t^{N_{\alpha,p}-3} |\omega_s(t)|^p dt   }$. Note that constant $C$ does not depend on $s$ and $S$.

\hfill $\Box$

Applying Lemma \ref{radialgenerlem1} enables us to prove the following pointwise estimate.

\begin{lem} \label{radialgenerlem2}
Let $n> 1-\alpha+\frac{\alpha}{p}$, $F\in C^1(\mathbb{R})$ and $u$ be a nonconstant radial stable weak solution of (\ref{gener}). Then, for fixed $\gamma>1$  there exists  $C_\gamma>0$ such that 
\begin{equation}\label{pointestimate}
|u(\gamma  r)-u(r)|\ge \frac{C_\gamma^{ \frac{p+1}{p}  }}{C^p} \ r^{N_p(\alpha)}   ,   \ \ \forall r\ge 1.
\end{equation}
The constant $C$ is the same as the constant in Lemma \ref{radialgenerlem1} which is independent of $\gamma$ and 
$$C_\gamma:=   \left\{
                      \begin{array}{ll}
                       ( 1+\frac{\alpha}{p}) \ln \gamma, & \hbox{if $N_{\alpha,p}=p+2$;} \\
                       \frac{p+1}{p+2-N_{\alpha,p}} \left(\gamma^{  \frac{ (p+2-N_{\alpha,p})(\alpha+p)   }{p(p+1)}  } -1\right), & \hbox{if $N_{\alpha,p}\neq p+2$;} \\
                           \end{array}
                    \right.$$  
                    where $N_{\alpha,p}:=\frac{p(n+\alpha)}{p+\alpha}$.
\end{lem}

\noindent\textbf{Proof:} Fix $\gamma>1$. By Lemma \ref{radialgenerlem1} for $1\le s\le \gamma s:=S$, we have 
\begin{eqnarray}\label{radiallemest1}
\int_{s}^{\gamma s} \frac{dt}{t^{N_{\alpha,p}-1}|\omega_s(t)|^p  } \le C s^{-2\sqrt{\frac{N_{\alpha,p}-1} {p-1}}} ,  
\end{eqnarray}
where $C$ does not depend on $\gamma$ and $s$.
On the other hand, by the same idea as in \cite{cc,ces,sv} we see $\omega_s$ does not change sign in $(0,\infty)$. So, 
\begin{eqnarray}\label{radiallemest2}
\int_{s}^{\gamma s} |\omega_s(t)| dt = |\omega(\gamma s)-\omega(s)|.
\end{eqnarray}
Now, apply the H\"{o}lder's inequality to get 
\begin{eqnarray*}
\int_{s}^{\gamma s} \frac{dt}{ t^{\frac{N_{\alpha,p}-1}{p+1}}  }&\le& \left(\int_{s}^{\gamma s}        \frac{dt}{t^{N_{\alpha,p}-1}|\omega_s(t)|^p  } \right)^{\frac{1}{p+1}}   \left(      \int_{s}^{\gamma s} |\omega_s(t)| dt       \right)^{\frac{p}{p+1}}\\
&\le& C^{\frac{1}{p+1}} s^{-\frac{2}{p+1}\sqrt{\frac{N_{\alpha,p}-1} {p-1}}} |\omega(\gamma s)-\omega(s)|^{\frac{p}{p+1}},
\end{eqnarray*}
 where in the last inequality we have used (\ref{radiallemest1}) and (\ref{radiallemest2}). On the other hand, by a direct calculation, one can see that  L.H.S. of the above inequality is
\begin{eqnarray*}
\int_{s}^{\gamma s} \frac{dt}{ t^{\frac{N_{\alpha,p}-1}{p+1}}  } = \hat C_\gamma\  s^{  \frac{p+2-N_{\alpha,p}}{p+1}   },
\end{eqnarray*}
where 
$$\hat C_\gamma:=   \left\{
                      \begin{array}{ll}
                        \ln \gamma, & \hbox{if $N_{\alpha,p}=p+2$;} \\
                       \frac{p+1}{p+2-N_{\alpha,p}} (\gamma^{  \frac{ p+2-N_{\alpha,p}   }{p+1}  } -1), & \hbox{if $N_{\alpha,p}\neq p+2$.} \\
                           \end{array}
                    \right.$$  
Therefore, 
\begin{eqnarray} \label{pointestimate2}
|\omega(\gamma s)-\omega(s)| \ge \frac{ \hat C_{\gamma}^{\frac{p+1}{p}}   }{ C^p}  s^{ \frac{1}{p}\left( p+2-N_{\alpha,p}  +2 \sqrt{\frac{N_{\alpha,p}-1} {p-1}} \right)}.
\end{eqnarray}
Since $\omega (s)=u(r)$ for $s=r^{1+\frac{\alpha}{p}}$ and $r>0$,  replace $\gamma $ by $\gamma^{ 1+\frac{\alpha}{p}   }$ in (\ref{pointestimate2}) to get 
\begin{eqnarray*} 
|u(\gamma r) - u(r)| \ge \frac{ \hat C_{\gamma^{   1+\frac{\alpha}{p} }   }^{\frac{p+1}{p}}   }{ C^p} \  r^{ \frac{1}{p}   \left(   1+\frac{\alpha}{p}   \right)       \left( p+2-N_{\alpha,p}  +2 \sqrt{\frac{N_{\alpha,p}-1} {p-1}} \right)},  \ \ \ \ \forall r>1.
\end{eqnarray*}

\hfill $\Box$ 

Now, we are in the position to prove Theorem \ref{radialgener2}.
\\
\\
\noindent\textbf{Proof of Theorem \ref{radialgener2}:} This is a consequence of Lemma \ref{radialgenerlem2}. Fix $1<\hat{\gamma} <\infty$. Since $u$ is bounded, R.H.S. of (\ref {pointestimate}) must be a bounded function of $r$. So, the exponent must be nonnegative, i.e.,  $N_p(\alpha) \le 0$.   Note that if the exponent is zero, then by the same idea as in the proof of Theorem \ref{radialgenerest} we get a lower bound of the $\log$-form. So, the exponent must be negative, and a straightforward calculation shows this is possible if and only if $n>p+\frac{4(p+\alpha)}{p-1}$. 

To get the desired pointwise estimate, we apply Lemma \ref{radialgenerlem2} again. Note that if $n>p+\frac{4(p+\alpha)}{p-1}$, then for sure
 $n>\frac{(p+2)(p+\alpha)}{p-1}-\alpha$, i.e. $N_{\alpha,p}>p+2$.  So, for any  $\gamma >1 $ from (\ref{pointestimate}), we conclude
\begin{eqnarray*} 
|u(\gamma r) - u(r)| \ge \frac{ C'_{\gamma }   } { C^p}  \ r^{ N_p(\alpha)},  \ \ \ \ \forall r>1,
\end{eqnarray*}
where $C'_\gamma= \left(\frac{p+1}{p+2-N_{\alpha,p}} (\gamma^{( \frac{ p+2-N_{\alpha,p}) (\alpha+p) }{(p+1)p}  } -1)\right)^{\frac{p+1}{p}}$.  Now, just take a limit of both sides of the above inequality when $\gamma\to \infty$. Note that $\lim_{\gamma\to \infty} C'_{\gamma}=\left(\frac{p+1}{N_{\alpha,p}-p-2}\right)^{\frac{p+1}{p}}<\infty$.   

\hfill $\Box$

\noindent\textbf{Proof of Theorem \ref{radialgenerest}:} If $n>p+\frac{4(p+\alpha)}{p-1}$, then a straightforward calculation shows $N_p(\alpha) < 0$.  Assume $\lim_{r\to \infty} u(r)=0$, otherwise (\ref {radialestimate}) holds trivially. Also, $u$ is monotone, because $u_r= (1+\frac{\alpha}{p}) r^{\frac{\alpha}{p}} \omega_s(s)$ and $\omega $ is monotone. Considering these facts and applying Lemma \ref{radialgenerlem2}, we get the following
\begin{eqnarray*}
|u(r)|&=& \sum_{k=1}^{\infty} |u(2^k r)-u(2^{k-1}r)| \ge C \sum_{k=1}^{\infty} (2^{k-1} r)^{ N_p(\alpha) }\\
&= & C \  r^{  N_p(\alpha) } \sum _{k=1}^{\infty} 2^{ (k-1)     N_p(\alpha) }.
\end{eqnarray*}
Since the exponent is negative, the series is convergent and (\ref {radialestimate}) holds. 

Now, let $1-\alpha+\frac{\alpha}{p}< n\le p+\frac{4(p+\alpha)}{p-1} $. Take $1 \le r_1<2$ such that $r=2^{m-1}r_1$. For $r\ge 1$, we have
\begin{eqnarray*}
|u(r)|&\ge & |u(r)-u(r_1)|-|u(r_1)| = \sum_{k=1}^{m-1}  |u(2^k r_1)-u(2^{k-1}r_1)| - |u(r_1)|\\
&\ge& C \sum_{k=1}^{m-1}  (2^{k-1} r_1)^{ N_p(\alpha)}  - |u(r_1)|.
\end{eqnarray*}
There are two cases.    If $n= p+\frac{4(p+\alpha)}{p-1}$, simple calculations show that $N_p(\alpha)=0$. So,  the latter inequality can be simplified to
\begin{eqnarray*}
|u(r)|&\ge & C(m-1)  - |u(r_1)| = C \frac{\ln r - \ln r_1}{\ln 2}  - |u(r_1)|. 
\end{eqnarray*}
Since $1\le r_1<2$ and $u$ is continuous function,  (\ref {radialestimate}) holds  for large enough $r$.  

Otherwise, for dimensions $1-\alpha+\frac{\alpha}{p}< n< p+\frac{4(p+\alpha)}{p-1} $, we have
\begin{eqnarray*}
|u(r)| \ge
 C \left( \frac{      r^{ N_p(\alpha) } -  r_1^{ N_p(\alpha)  }      }{       2^{  N_p(\alpha) }   -1  } \right) -|u(r_1)|.
\end{eqnarray*}
Since  $  N_p(\alpha)>0$, we get the desired result for large enough $r$.

\hfill $\Box$

  Note that regularity of extremal solutions which is closely related to Liouville theorems through blow-up analysis and rescaling techniques for both semilinear and quasilinear equations has been established by Cabr\'e et al. in \cite{cc1,ccs}.

%%%%%%%%%%%%%%%%%%%%%%%%%%%%%%%%%%%%%%%%%%%%%%%

\subsection{Proof of Theorem \ref{stablenon}.}\label{2}

 At first, in a couple of lemmata, we  prove some major estimates for stable weak solutions  (not necessarily radial) of (\ref{gener}) with three different nonlinearities $F(u)=e^u$, $u^q$ where $q>p-1$ and $-u^{q}$ where $q<0$,  known as the Gelfand, the Lane-Emden and the negative exponent nonlinearities, respectively.  Note that the negative exponent nonlinearity for $q=-2$ is called the MEMS nonlinearity, see \cite{gg,egg2}. Then, applying appropriate test functions leads us to Liouville theorems. The following lemmata are adaptations of the same type estimates given by Farina in \cite{f1,f2} and Esposito-Ghoussoub-Guo in \cite{egg1,egg2}. Similar results can be found in \cite{ces,df,df1,e1,e2,cf}. 

For power nonlinearities, we have:
\begin{lem}  \label{stpower1}
Let $\Omega\subset \mathbb{R}^n$ and $u\in C^{1,\gamma}(\Omega)$ be positive stable weak solution of (\ref{gener}) with $F(u)=\text{sign} (q) u^q$. For one of the following exponents and parameters:
\\
(i) either $q>p-1$  and  $ 1\le t<-1+2\frac{q+\sqrt{q(q-p+1)}}{p-1}$,\\
(ii) or $q<0$ and $ 1\le -t<1+2\frac{-q+\sqrt{q(q-p+1)}}{p-1}$. Then, we have
\begin{eqnarray} \label{estpower1}
\int_{\Omega}\left( |\nabla u|^p u^{t-1}+ f(x)  u^{t+q}\right)\phi^{pm} &\le&C \int_{\Omega} f(x)^{-\frac{t+p-1}{q-p+1}  }| \nabla \phi|^{\frac{t+q}{q-p+1}p},
\\
 \label{estpower2}\int_{\Omega} \left( |\nabla u|^p u^{t-1}+f(x) u^{t+q}\right)\phi^{2m} &\le&  C \int_{\Omega} f(x)^{-\frac{t+1}{q-1}   }(|\nabla u|^{p-2} | \nabla \phi|^2)^{\frac{t+q}{q-1}},
\end{eqnarray}
for all $\phi\in C_c^1(\Omega)$ with $0\le\phi\le 1$ and large enough $m$. The constant $C$ does not depend on $\Omega$ and $u$.
\end{lem}

\begin{remark}\label{remarkM}
If we assume that $f\ge M>0$ where $M$ is a constant, then $f(x)$ in the right hand sides of (\ref{estpower1}) and (\ref{estpower2}) can be replaced by  $M^{-\frac{t+p-1}{q-p+1}  }$ and $M^{-\frac{t+1}{q-1}   }$, respectively. Therefore, strictly positive weight $f(x)$ does not have any impact on these decay estimates. 
\end{remark}

\noindent\textbf{Proof:} We first prove (\ref{estpower1}), then by the same idea we prove (\ref{estpower2}). For any $C^{1,\gamma}(\mathbb{R}^n)$ stable solution of (\ref{gener}) with $F(u)=\text{sign}(q)u^q$ and $\phi\in C_c^1(\mathbb{R}^n)$, we have the followings:
\begin{eqnarray}\label{stab}
|q| \int_{\Omega}  f(x) u^{q-1} \phi^2 &\le& (p-1) \int_{\Omega} |\nabla u |^{p-2} |\nabla \phi|^2,\\
\label{pde}\text{sign}(q) \int f(x) u^q \phi &=& \int |\nabla u|^{p-2} \nabla u\cdot\nabla \phi.
\end{eqnarray} 

Test (\ref{pde}) on $u^t \phi^p$, an appropriate $t\in \mathbb{R}$ will be chosen later,  to get 
\begin{eqnarray*}
\text{sign}(q) \int_{\Omega} f(x) u^{t+q} \phi^p &=& \int |\nabla u|^{p-2} \nabla u\cdot\nabla \left(u^t \phi^p\right)\\
& = &t \int_{\Omega} |\nabla u|^p u^{t-1} \phi^p + p \int |\nabla u |^{p-2} u^t \nabla u\cdot\nabla \phi \phi^{p-1}.
\end{eqnarray*} 
Therefore, applying the Young's inequality\footnote{For $a,b,\epsilon>0$ and $1<\alpha,\beta<\infty$ we have $ab\le \epsilon a^\alpha+C(\epsilon) b^\beta$, where $C(\epsilon)=(\epsilon\alpha)^{-\beta/\alpha} \beta^{-1}$ and $1/\alpha+1/\beta=1$.  For $\alpha=\beta=2$ this is called the Cauchy's inequality.} with exponents $p$ and $\frac{p}{p-1}$ to $\left(  |\nabla u|^{p-1}  u^{\frac{p-1}{p}(t-1)} \phi^{p-1} \right)\left( u^{\frac{t+p-1}{p}} |\nabla \phi|  \right)$, we have 
\begin{eqnarray}\label{pde10}
(|t|-\epsilon) \int |\nabla u|^p u^{t-1} \phi^p \le C_{\epsilon,p} \int  u^{t+p-1}|\nabla \phi|^p + \int f(x)  u^{t+q} \phi^p.
 \end{eqnarray} 
Now, test (\ref{stab}) on $u^{\frac{t+1}{2}}\phi^{\frac{p}{2}}$ to obtain
\begin{eqnarray*}
\frac{|q|}{p-1}  \int f(x) u^{t+q} \phi^p &\le& \frac{(t+1)^2}{4} \int |\nabla u|^p  u^{t-1} \phi^p +\frac{p^2}{4} \int |\nabla u|^{p-2}u^{t+1} \phi^{p-2} |\nabla\phi|^2  \\
&&+ \frac{(t+1)p}{2} \int |\nabla u|^{p-2} u^t  \nabla u \cdot \nabla \phi \phi^{p-1} \\
&\le& \left(   \frac{(t+1)^2}{4} +2\epsilon  \right) \int |\nabla u|^p u^{t-1} \phi^p + (C'_{\epsilon,t,p}+C''_{\epsilon,t,p} ) \int   u^{t+p-1}|\nabla \phi |^p,
 \end{eqnarray*} 
in the last inequality we have used the Young's inequality twice with exponents $p$ and $\frac{p}{p-1}$ and also with $\frac{p}{2}$ and $\frac{p}{p-2}$. Combine this inequality and (\ref{pde10}) to see
\begin{eqnarray} \label{majorest}
\left(    \frac{|q|}{p-1} -\frac{ \frac{(t+1)^2}{4} +2\epsilon  }{|t|-\epsilon}    \right)
\int f(x) u^{t+q} \phi^p \le \left(   \frac{     \frac{(t+1)^2}{4} +2\epsilon  }{|t|-\epsilon} C_{\epsilon,p}  +C'_{\epsilon,t,p}+C''_{\epsilon,t,p}\right) \int  u^{t+p-1}|\nabla \phi|^p.
 \end{eqnarray}
For an  appropriate choice of $t$, given in the assumption, we see that the coefficient in L.H.S. is positive for $\epsilon$ small enough. Therefore, replacing $\phi $ with $\phi^m$ for large enough $m$ and applying the H\"{o}lder's  inequality  with exponents $\frac{t+q}{t+p-1}$ and $\frac{t+q}{q-p+1}$ we obtain 
\begin{eqnarray} \label{majorest1}
\int f(x) u^{t+q} \phi^{pm} \le D_{\epsilon,t,m,p} \int f(x)^{-\frac{t+p-1}{q-p+1}}  |\nabla\phi |^{\frac{t+q}{q-p+1}p}. 
 \end{eqnarray}
Note that both exponents are greater than 1 for $t$ given in (i) and (ii).

On the other hand, combining (\ref{pde10}) and (\ref{majorest}) gives us
\begin{eqnarray*}
\int |\nabla u|^p u^{t-1} \phi^p \le D'_{\epsilon,t,p} \int  u^{t+p-1}|\nabla \phi|^p.
 \end{eqnarray*}
Similarly, replace $\phi $ by $\phi^m$ and apply H\"{o}lder's  inequality with exponents $\frac{t+q}{t+p-1}$ and $\frac{t+q}{q-p+1}$ to get  
\begin{eqnarray*}
\int |\nabla u|^p u^{t-1} \phi^{mp} \le D''_{\epsilon,t,p,m} \int  f(x)^{ -\frac{t+p-1}{q-p+1}}   |\nabla \phi|^{  \frac{t+q}{q-p+1}p}.
 \end{eqnarray*}
This inequality and (\ref{majorest1}) finish the proof of (\ref{estpower1}). 

Proof of (\ref{estpower2}) is quite similar. Here is the sketch of proof. First, test (\ref{pde}) on $u^t\phi^2$ to arrive at 
\begin{eqnarray*}
\text{sign}(q) \int f(x) u^{t+q} \phi^2 &=&  t \int |\nabla u|^p u^{t-1} \phi^2 + 2 \int |\nabla u |^{p-2} u^t  \nabla u\cdot\nabla \phi   \ \phi.
\end{eqnarray*} 
Then, applying Cauchy's inequality to $\left(|\nabla u|^{\frac{p}{2}} u^{\frac{t-1}{2}}\phi\right)\left(|\nabla u|^{\frac{p-2}{2}} u^{\frac{t+1}{2}} |\nabla \phi|\right)$, we have 
\begin{eqnarray}\label{pde2}
(|t|-\epsilon) \int |\nabla u|^p u^{t-1} \phi^2 \le \hat{C}_{\epsilon,p} \int  |\nabla u|^{p-2} u^{t+1} |\nabla\phi|^2 + \int f(x)  u^{t+q} \phi^2.
 \end{eqnarray} 
Now, test (\ref{stab}) on $u^{\frac{t+1}{2}} \phi$ and apply the Cauchy's inequality to get
\begin{eqnarray*}
\frac{|q|}{p-1}  \int f(x) u^{t+q} \phi^2 &\le& \frac{(t+1)^2}{4} \int |\nabla u|^p  u^{t-1} \phi^2 + \int u^{t+1}  |\nabla u|^{p-2}|\nabla\phi|^2 \\
&& + (t+1) \int |\nabla u|^{p-2} u^t  \nabla u \cdot \nabla \phi \phi \\
&\le& \left(   \frac{(t+1)^2}{4} +\epsilon  \right) \int |\nabla u|^p u^{t-1} \phi^2 + (1+\hat{C}'_{\epsilon,t,p} ) \int  |\nabla u|^{p-2}u^{t+1}|\nabla \phi |^2.
 \end{eqnarray*} 
 Combine this inequality and (\ref{pde2}) to see
\begin{eqnarray*} 
\left(    \frac{|q|}{p-1} -\frac{ \frac{(t+1)^2}{4} +\epsilon  }{|t|-\epsilon}    \right)
\int f(x) u^{t+q} \phi^2 \le \left( 1+\hat{C}'_{\epsilon,t,p} +\frac{\hat{C}_{\epsilon,p}}{|t|-\epsilon} \left(\frac{(t+1)^2}{4}+\epsilon\right) \right) \int |\nabla u|^{p-2} u^{t+1}|\nabla \phi|^2.
 \end{eqnarray*}
Replace $\phi$ with $\phi^m$ and apply the H\"{o}lder's inequality to get the desired result.

\hfill $\Box$

For the exponential nonlinearity that is Gelfand nonlinearity, we have the following estimate.

\begin{lem} \label{stgel1}
Let $\Omega\subset \mathbb{R}^n$ and $u\in C^{1,\gamma}(\Omega)$ be stable weak solution of (\ref{gener}) with $F(u)=e^u$. For any $0<t<\frac{2}{p-1}$, we have 
\begin{eqnarray} \label{estgel1}
\int_{\Omega}  f(x) e^{(2t+1)u}\phi^{pm} &\le & C \int_{\Omega}  f(x)^{-2t}| \nabla \phi|^{p(2t+1)},
\\
\label{estgel2} \int_{\Omega}  f(x) e^{(2t+1)u}\phi^{2m} &\le & C \int_{\Omega}  f(x)^{-2t} (|\nabla u|^{p-2}| \nabla \phi|^2)^{(2t+1)},
\end{eqnarray}
for all $\phi\in C_c^1(\Omega)$ with $0\le\phi\le 1$ and large enough $m$. The constant $C$ does not depend on $\Omega$ and $u$.
\end{lem}

\noindent\textbf{Proof:} The idea of the proof is the same as Lemma \ref{stpower1}.  We just prove (\ref{estgel1}), then by the same idea one can prove (\ref{estgel2}). For any stable solution of (\ref{gener}) with $F(u)=e^u$ and $\phi\in C_c^1(\mathbb{R}^n)$, we have the following:
\begin{eqnarray}\label{stab1}
 \int_\Omega  f(x) e^u \phi^2 &\le& (p-1) \int_\Omega |\nabla u |^{p-2} |\nabla \phi|^2,\\
\label{pde1} \int_\Omega f(x) e^u \phi &=& \int_\Omega |\nabla u|^{p-2} \nabla u\cdot\nabla \phi.
\end{eqnarray} 

Test (\ref{pde1}) on $e^{2tu} \phi^p$, an appropriate $t\in \mathbb{R}^+$ will be chosen later,  to get 
\begin{eqnarray*}
 \int f(x) e^{(2t+1)u} \phi^p &=& \int |\nabla u|^{p-2} \nabla u\cdot\nabla \left( e^{2tu} \phi^p\right)\\
&=& 2t \int |\nabla u|^p e^{2tu} \phi^p + p \int |\nabla u |^{p-2} e^{2tu} \nabla u\cdot\nabla \phi \ \phi^{p-1}.
\end{eqnarray*} 
Apply the Young's inequality with exponents $p$ and $\frac{p}{p-1}$ to $\left(  |\nabla u|^{p-1}  e^{\frac{2t(p-1)}{p}u} \phi^{p-1} \right)\left( e^{\frac{2t}{p}u} |\nabla \phi|  \right)$ to obtain
\begin{eqnarray}\label{pde11}
(2t-\epsilon) \int |\nabla u|^p e^{2tu} \phi^p \le \tilde{C}_{\epsilon,p} \int  e^{2tu} |\nabla \phi|^p + \int f(x)  e^{(2t+1)u} \phi^p.
 \end{eqnarray} 
Now, test (\ref{stab1}) on $e^{tu} \phi^{\frac{p}{2}}$ to have
\begin{eqnarray*}
\frac{1}{p-1}  \int f(x) e^{(2t+1)u} \phi^p &\le& t^2 \int |\nabla u|^p  e^{2tu} \phi^p +\frac{p^2}{4} \int |\nabla u|^{p-2} e^{2tu} \phi^{p-2} |\nabla\phi|^2  \\
&&+ tp \int |\nabla u|^{p-2} e^{2tu}  \nabla u \cdot \nabla \phi \phi^{p-1} \\
&\le& \left(  t^2 +2\epsilon  \right) \int |\nabla u|^p e^{2tu} \phi^p + (\tilde{C}'_{\epsilon,t,p}+\tilde{C}''_{\epsilon,t,p} ) \int   e^{2tu} |\nabla \phi |^p,
 \end{eqnarray*} 
in the last inequality we have used the Young's inequality twice with exponents $p$ and $\frac{p}{p-1}$ and also with $\frac{p}{2}$ and $\frac{p}{p-2}$. Combine this inequality and (\ref{pde11}) to see
\begin{eqnarray} \label{majorest11}
\left(    \frac{1}{p-1} -\frac{ t^2 +2\epsilon  }{2t-\epsilon}    \right)
\int f(x) e^{(2t+1)u} \phi^p \le \left(     \frac{t^2+2\epsilon}{2t-\epsilon} \tilde{C}_{\epsilon,p} +\tilde{C}'_{\epsilon,t,p}+\tilde{C}''_{\epsilon,t,p}\right) \int  e^{2tu}|\nabla \phi|^p.
 \end{eqnarray}
For $\epsilon$ small enough,  choosing $0<t<\frac{2}{p-1}$ we see that the coefficient in the L.H.S. is positive. 

Now, replacing $\phi $ with $\phi^m$ for large enough $m$ and applying the H\"{o}lder's  inequality  with exponents $2t+1$ and $\frac{2t+1}{2t}$ we obtain 
\begin{eqnarray*} 
\int f(x) e^{(2t+1)u} \phi^{pm} \le \tilde{D}_{\epsilon,t,p,m} \int f(x)^{-2t}  |\nabla\phi |^{p(2t+1)}. 
 \end{eqnarray*}

\hfill $\Box$

To prove the theorem we just pick an appropriate test function.
\\
\\
\noindent\textbf{Proof of Theorem \ref{stablenon}:} We only prove the results for the exponential nonlinearity. Let $\zeta_R \in C^1_c(\mathbb{R}^n)$ such that $0\le \zeta_R\le 1$ be given by $$\zeta_R(x)=\left\{
                      \begin{array}{ll}
                        1, & \hbox{if $|x|<R$;} \\
                        0, & \hbox{if $|x|>2R$;} \\
                           \end{array}
                    \right.$$                             
and $||\nabla\zeta_R||_{\infty}\le\frac{\hat C}{R}$.  Test (\ref{estgel1}) on $\zeta_R$ to get 
\begin{eqnarray*}
\int_{B_R} |x|^\alpha e^{(2t+1)u}  \le C \int_{B_{2R}} |x|^{-2t\alpha}  |\nabla\zeta_R |^{p(2t+1)}= C_{n,\alpha,t,p} R^{n-2t\alpha-p(2t+1)},
 \end{eqnarray*}
where $0<t<\frac{2}{p-1}$. If $n<2t(\alpha+p)+p$ for $0<t<\frac{2}{p-1}$, i.e. $n<\frac{4(p+\alpha)}{p-1}+p$, by sending $R\to\infty$ we get the contradiction.

\hfill $\Box$

\subsection{Proofs of Theorem \ref{radialcompact} and \ref{poho1}.}\label{3}

 The crucial tool to prove these theorems is the following Harnack's inequality, see \cite{gt,s1}. 

\begin{lem}\label{har} (Harnack's inequality) Assume $\Omega\subset \mathbb{R}^n$ and $n>p$. Let $w$ be a nonnegative weak subsolution of
\begin{equation*}
-\Delta_p w= a(x) w^{p-1} \ \ \text{in } \ \Omega,
\end{equation*}
where $a(x)\in L^q(\Omega)$ for $q\in (\frac{n}{p},\frac{n}{p-1})$. Then,  for any $R$ such that $B_{2R}\subset \Omega$, there exists $C_H$ such that
\begin{equation*}
||w||_{L^{\infty}(\Omega)} \le C_H R^{-\frac{n}{\beta}} || w||_{L^{\beta}(\Omega)},
\end{equation*}
where $\beta>1$ and $C_H$ may depend on $p,n, \beta$ and $R^{p-\frac{n}{q}}||a||_{L^{q}}$.
\end{lem}

\noindent\textbf{Proof of Theorem \ref{radialcompact}:} 
In the light of Proposition \ref{change}, we can assume $\alpha=0$. Let $\Sigma\subset B_{R_0}$ for sufficiently large enough $R_0$.

(i) Let $F(u)=e^u$. Apply Lemma \ref{stgel1} with $\Omega=\mathbb{R}^n\setminus \overline{B_{R_0}}$ and the following test function $\xi_R\in C^1_c(\mathbb{R}^n\setminus \Sigma)$ for $R>R_0+3$;
 $$\xi_R(x)=\left\{
                      \begin{array}{ll}
                      0, & \hbox{if $|x|<R_0+1$;} \\
                        1, & \hbox{if $R_0+2<|x|<R$;} \\
                        0, & \hbox{if $|x|>2R$;} 
                                                                       \end{array}
                    \right.$$
which satisfies $0\le\xi_R\le1$, $||\nabla\xi_R||_{L^{\infty}(B_{2R}\setminus B_{R})}<\frac{C}{R}$ and $||\nabla\xi_R||_{L^{\infty}(B_{R_0+2}\setminus B_{R_0+1})}<C_{R_0}$. Therefore, for $R>R_0+3$ and $0<t<\frac{2}{p-1}$, we get 
\begin{eqnarray}\label{mainest}
\int_{B_R\setminus B_{R_0+2}}  e^{(2t+1)u} \le C_{R_0}+C  R^{n-p(2t+1)}.
\end{eqnarray}
 Since $p<n<\frac{4p}{p-1}+p$, we can choose $t_1:=\frac{n-p}{2p}$ in (\ref{mainest}).   By sending $R$ to infinity, we see 
\begin{eqnarray*}
\int_{\mathbb{R}^n\setminus B_{R_0+2}}  e^{\frac{n}{p}u} < \infty.
\end{eqnarray*}
So, for a given $\delta>0$ and large enough $R_1>R_0+3$, we obtain 
\begin{eqnarray}\label{mainest1}
\int_{\mathbb{R}^n\setminus B_{R_1}}  e^{\frac{n}{p}u} \le \delta^{\frac{n}{p}}.
\end{eqnarray}
Now, take $B_{2R}(y)\subset \{x; |x|>R_1\}\subset \{x; |x|>R_0\}$ and $|y|=4R$  for $R>R_0$. Then, we have  $B_{2R}(y)\subset \{x; \ 2R<|x|<6R\}$. Using standard test functions of the form $\phi_R(x)=\zeta_R({|x-y|})$ for $\zeta_R\in C^1_c(\mathbb{R}^n\setminus B_{R_0})$ which satisfies  $0\le\zeta_R\le1$, $||\nabla\zeta_R||_{\infty}<\frac{C}{R}$ and
 $$\zeta_R(x)=\left\{
                      \begin{array}{ll}
                        1, & \hbox{if $|x|<R$;} \\
                        0, & \hbox{if $|x|>2R$,} 
                       \end{array}
                    \right.$$
in (\ref{estgel1}), we get  
\begin{eqnarray}\label{mainest2}
\int_{B_R(y)} e^{(2t+1)u} \le C R^{n-p(2t+1)},
\end{eqnarray}
for $0<t<\frac{2}{p-1}$. The positive constant $C$ is independent of $R$.

To apply the Harnack's inequality, set $w:=e^u$  and observe that $w$ is a positive subsolution of the following equation $-\Delta_p w = e^u w^{p-1}$ in $B_{2R}(y)$.  Therefore, in the light of Lemma \ref{har} for $\beta=\frac{n}{p}>1$, we observe
\begin{eqnarray*}
\sup_{B_{2R}(y)} e^u &\le &  C_H R^{-\frac{n}{\beta}} \left(\int_{B_{2R}} e^{\frac{n}{p}u}  \right)^{\frac{1}{\beta}}\\
&\le& C_H R^{-\frac{n}{\beta}} \delta,
\end{eqnarray*}
in the last inequality we have used (\ref{mainest1}). Therefore,  
\begin{eqnarray}\label{mainest3}
\sup_{B_{2R}(y)} e^u \le C_H \delta R^{-p},
\end{eqnarray}
where $C_H$ just depends on $n,p$ and $R^{\epsilon} || e^u  ||_{L^{\frac {n}{p-\epsilon}}(B_{2R}(y))}$. Set $t_2:=\frac{n}{2(p-\epsilon)}-\frac{1}{2}$ in (\ref{mainest2}) for $\epsilon $ small enough, to get
\begin{eqnarray*}
R^{\epsilon} ||   e^u  ||_{L^{\frac {n}{p-\epsilon}}(B_{2R}(y))} &\le& 
R^{\epsilon}  \left (   \int_{B_{2R}(y)}  e^{\frac{n}{p-\epsilon}u} \right)^{\frac{p-\epsilon}{n}}\le C R^\epsilon R^{\frac{p-\epsilon}{n} \left(n-\frac{n}{p-\epsilon}p  \right)}=C,
\end{eqnarray*}
where $C$ only depends on $n$ and $p$.  Since $|y|=4R$,  from (\ref{mainest3}), we see $ |y|^p e^{u(y)} \le \delta$ for a given $\delta$ and large enough $|y|$. Hence, 
\begin{eqnarray}\label{decay}
\lim_{|y|\to \infty} |y|^p e^{u(y)}=0.
\end{eqnarray}
Note that we have not used the assumption that $u$ is radial so far. 

For a radial weak solution $u$ of (\ref{gener}),  the definition of radial $p$-Laplacian operator and decay estimate (\ref{decay}), imply there exists $R_2$ and $0<k<n-p$ such that 
\begin{eqnarray*}
-r^{1-n}(r^{n-1} |u_r|^{p-2}u_r )_r\le k r^{-p},  \ \ \ \forall r>R_2.
\end{eqnarray*}
By integration,  we get 
\begin{eqnarray*}
 |u_r|^{p-2}u_r \ge -\frac{k}{n-p}r^{1-p} +\frac{C(n)}{r^{n-1}},  \ \ \ \forall r>R_2.
\end{eqnarray*}
So, for large enough $R_3$ we have 
\begin{eqnarray*}
  |u_r|^{p-2}u_r \ge - r^{1-p} ,  \ \ \ \forall r>R_3,
\end{eqnarray*}
and by integrating this, we get
\begin{eqnarray*}
 r^{p} e^{u(r) }\ge C \ r^{p-1},   \ \ \ \forall r>R_3.
\end{eqnarray*}
Since $p\ge 2$, this is in contradiction with (\ref{decay}). Hence, there is no radial stable outside a compact set solution for (\ref{gener}) with $\alpha=0$ in dimensions $p<n<\frac{4p}{p-1}+p$. Now, apply Proposition \ref{change} to see there is no such solution for (\ref{gener}) with any $\alpha>-p$ in dimensions $p<n<\frac{4(p+\alpha)}{p-1}+p$.

(ii)  The case $F(u)=-u^{q}$ for $q>p-1$ and $\alpha=0 $ has been done in \cite{dfsv}.

(iii) Let $F(u)=u^{q}$ for $q<0$ and $\alpha=0$. By a similar argument as in (i), i.e., applying Lemma \ref{stpower1} with the same test function $\xi_R$, we obtain 
\begin{eqnarray}\label{mainestm}
\int_{B_R\setminus B_{R_0+2}}  u^{t+q} \le C_{R_0}+C  R^{n-\frac{t+q}{q-p+1}p},
\end{eqnarray}
for $1\le -t<1+2\frac{-q+\sqrt{q(q-p+1)}}{p-1}$.  Since $\frac{p(q-1)}{q-p+1}<n<\frac{q-1}{q-p+1}p+2p\frac{q-\sqrt{q(q-p+1)}}{(p-1)(q-p+1)}$, we can take $ t_1:=\frac{n}{p}(q-p+1)-q$ in (\ref{mainestm}).   By sending $R$ to infinity, we see 
\begin{eqnarray*}
\int_{\mathbb{R}^n\setminus B_{R_0+2}}  u^{\frac{n}{p}(q-p+1)} < \infty.
\end{eqnarray*}
So, for a given $\delta>0$ and large enough $R_1>R_0+3$, we have 
\begin{eqnarray}\label{mainestm1}
\int_{\mathbb{R}^n\setminus B_{R_1}} u^{\frac{n}{p}(q-p+1)} \le \delta^{\frac{n}{p}(q-p+1)}.
\end{eqnarray}

On the other hand, applying Lemma \ref{stpower1} with the same test function, $\zeta_R$, as in (i), we get 
\begin{eqnarray}\label{mainestm2}
\int_{B_R(y)}  u^{t+q} \le C  R^{n-\frac{t+q}{q-p+1}p},
\end{eqnarray}
for $1\le -t<1+2\frac{-q+\sqrt{q(q-p+1)}}{p-1}$.

Now, define $w:=u^{-1}$ and observe by a straightforward calculation that 
$$\Delta_p w+u^{q-p+1} w^{p-1}= 2(p-1) u^{-2p+1} |\nabla u|^p.$$ Therefore, $w$ is a positive subsolution for $-\Delta_p w=a(x) w^{p-1}$, where $a(x)=u^{q-p+1}$. Apply the Harnack's inequality with  $\hat \beta:=-\frac{n}{p}(q-p+1)>1$, to get 
\begin{eqnarray*}
\sup_{B_{2R}(y)} w &\le &   {C}_H R^{-\frac{n}{\hat\beta}}   \left(\int_{B_{2R}} u^{\frac{n}{p}(q-p+1)} \right )^{\frac{1}{\hat\beta}}\\
&\le& C_H  \delta R^{-\frac{n}{\hat\beta}},
\end{eqnarray*}
where in the last inequality we have used (\ref{mainestm1}). Substitute the value of $\hat \beta$ to arrive at
\begin{eqnarray}\label{mainestm3}
\sup_{B_{2R}(y)} u^{-1} \le C_H \delta R^{\frac{p}{q-p+1}},
\end{eqnarray}
where $C_H$ just depends on $n,p$ and $R^{\epsilon} || u^{q-p+1}  ||_{L^{\frac {n}{p-\epsilon}}(B_{2R}(y))}$. Set $ t_2:=\frac{n}{p-\epsilon}(q-p+1)-q$ in (\ref{mainestm2}) for $\epsilon $ small enough, to get
\begin{eqnarray*}
R^{\epsilon} ||   u^{q-p+1}  ||_{L^{\frac {n}{p-\epsilon}}(B_{2R}(y))} &\le& 
R^{\epsilon}  \left (   \int_{B_{2R}(y)}  u^{\frac{n}{p-\epsilon}(q-p+1)} \right)^{\frac{p-\epsilon}{n}}\le C R^\epsilon R^{\frac{p-\epsilon}{n} \left(n-\frac{n}{p-\epsilon}p  \right)}=C,
\end{eqnarray*}
where $C$ only depends on $n$ and $p$.  Since $|y|=4R$,  (\ref{mainestm3}) proves the following decay estimate for not necessary radial solutions, 
\begin{eqnarray}\label{decaym}
\lim_{|y|\to \infty} |y|^{\frac{-p}{q-p+1}} {u^{-1}(y)}=0.
\end{eqnarray}

Note that for a radial solution of (\ref{gener}) we have $r^{n-1}|u_r|^{p-2} u_r$ is increasing and $u_r>0$. From the decay estimate (\ref{decaym}),  there exist $R_1$ and $k>0$ such that 
\begin{eqnarray*}
r^{1-n}(r^{n-1} |u_r|^{p-2}u_r )_r\le k r^{-\frac{pq}{q-p+1}},  \ \ \ \forall r>R_1.
\end{eqnarray*}
By integration,  we get 
\begin{eqnarray*}
 r^{n-1}u_r^{p-1} \le \frac{k}{n-\frac{pq}{q-p+1}}r^{n-\frac{pq}{q-p+1}} +C(n),  \ \ \ \forall r>R_1.
\end{eqnarray*}
Since $n>\frac{pq}{q-p+1}$, there exists positive constant $C$ independent of $r$ such that for large enough $R_2$, we have 
\begin{eqnarray*}
  u_r^{p-1} \le C r^{-\frac{(q+1)(p-1)}{q-p+1}} ,  \ \ \ \forall r>R_2,
\end{eqnarray*}
and again by integration, we get
\begin{eqnarray*}
  u(r) \le \hat C \ r^{-\frac{p}{q-p+1}},   \ \ \ \forall r>R_3,
\end{eqnarray*}
where $\hat C $ is a positive constant independent of $r$. This is in contradiction with (\ref{decaym}). Hence, there is no radial stable outside a compact set solution for (\ref{gener}) with $\alpha=0$ in dimensions $\frac{p(q-1)}{q-p+1}<n<\frac{q-1}{q-p+1}p+2p\frac{q-\sqrt{q(q-p+1)}}{(p-1)(q-p+1)}$. Now, apply Proposition \ref{change} to see there is no such solutions for (\ref{gener}) in the given dimension for any $\alpha+p>0$.

\hfill $\Box$ 

 \begin{remark}
We would like to emphasize that in the above theorem we have followed the ideas that are initiated by Dancer and Farina  \cite{df} for the semilinear Gelfand equation that is (\ref{gener}) with $\alpha=0$, $p=2$ and $F(u)=e^u$. Similar methods have been used by Farina in \cite{f1} for the  semilinear Lane-Emden equation that is (\ref{gener}) with $\alpha=0$,  $p=2$ and $F(u)=|u|^{q-1}u$. Recently, Wang and Ye in \cite{wy} applied these methods to the weighted semilinear Gelfand and Lane-Emden equations. For the negative exponent nonlinearity, we refer to Esposito-Ghoussoub-Guo \cite{egg2}.  The  quasilinear Lane-Emden equation that is (\ref{gener}) with $\alpha=0$, $p>2$ and $F(u)=u^q$ when $q>p-1$ has been studies by Damascelli et al. in \cite{dfsv}.

 \end{remark}

\noindent\textbf{Proof of Theorem \ref{poho1}.} The proof is an adaptation of the proof of Theorem 1.9 and Proposition 10.1 in \cite{dfsv}.
\\
Step 1: It's straightforward to observe the following Pohozaev type identity holds on any $\Omega \subset\mathbb{R}^n$. 
\begin{equation}\label{pohozaevidentity}
\frac{n+\alpha}{q+1} \int_{\Omega} |x|^\alpha u^{q+1} - \frac{n-p}{p} \int_{\Omega} |\nabla u|^p = \frac{1}{q+1} \int_{\partial \Omega} |x|^\alpha u^{q+1} x\cdot\nu + \int_{\partial\Omega} |\nabla u |^{p-2} x\cdot\nabla u\ u_\nu -\frac{1}{p} \int_{\partial\Omega} |\nabla u|^p x\cdot\nu .
\end{equation}
Step 2: Estimates $|\nabla u|\in L^p(\mathbb{R}^n)$ and $|x|^\alpha u^{q+1} \in L^1(\mathbb{R}^n)$ hold. To prove this, we use (\ref{estpower1}) with the following test function $\xi_R\in C^1_c(\mathbb{R}^n\setminus \Sigma)$ for $R>R_0+3$ and $\Sigma \subset B_{R_0}$;
 $$\xi_R(x):=\left\{
                      \begin{array}{ll}
                      0, & \hbox{if $|x|<R_0+1$;} \\
                        1, & \hbox{if $R_0+2<|x|<R$;} \\
                        0, & \hbox{if $|x|>2R$;} 
                                                                       \end{array}
                    \right.$$
which satisfies $0\le\xi_R\le1$, $||\nabla\xi_R||_{L^{\infty}(B_{2R}\setminus B_{R})}<\frac{C}{R}$ and $||\nabla\xi_R||_{L^{\infty}(B_{R_0+2}\setminus B_{R_0+1})}<C_{R_0}$. Therefore, 
$$\int_{R_0+2<|x|<R} (|\nabla u|^p u^{t-1} + |x|^\alpha u^{t+q})\le C_{R_0}+\hat C\ R^{n- \frac{p(t+q)}{q-p+1}   -\frac{t+p-1}{q-p+1} \alpha  }, $$ 
for all $1\le t<-1+2\frac{q+\sqrt{q(q-p+1)}}{p-1}$. Now,  set $t=1$ and send $R\to \infty$. Since $n< \frac{p(q+\alpha+1)}{q-p+1}$, we see
$\int_{\mathbb{R}^n} |\nabla u|^p<\infty$ and $\int_{\mathbb{R}^n}|x|^\alpha u^{q+1}<\infty$.
\\
\\
Step 3:  This equality holds $$(\frac{n-p}{p}-\frac{n+\alpha}{q+1})\int_{\mathbb{R}^n} |x|^\alpha u^{q+1}=0.$$  
Multiply (\ref{gener}) with $u\zeta_R$ for $\zeta_R\in C^1_c(\mathbb{R}^n)$ which satisfies  $0\le\zeta_R\le1$, $||\nabla\zeta_R||_{\infty}<\frac{C}{R}$ and
 $$\zeta_R(x):=\left\{
                      \begin{array}{ll}
                        1, & \hbox{if $|x|<R$;} \\
                        0, & \hbox{if $|x|>2R$.} 
                                                                       \end{array}
                    \right.$$
Then, integrate over $B_{2R}$ to get 
\begin{eqnarray} \label{identitytest}
\int_{B_{2R}} |x|^\alpha u^{q+1} \zeta_R - \int_{B_{2R}} |\nabla u|^p \zeta_R = \int_{B_{2R}} |\nabla u|^{p-2} \nabla \zeta_R\cdot \nabla u \ u.
\end{eqnarray}
By H\"{o}lder's inequality, we have the following upper bound for R.H.S. of (\ref{identitytest}),
\begin{eqnarray*}
|\int_{B_{2R}} |\nabla u|^{p-2} \nabla \zeta_R\cdot \nabla u \ u| & \le & R^{-1} \int_{B_{2R}} |\nabla u|^{p-1} ( |x|^{\frac{\alpha}{q+1}}u )\ |x|^{-\frac{\alpha}{q+1}}\\
&\le& R^{-1}\left(\int_{B_{2R}}  |\nabla u|^p \right)^{\frac{p-1}{p}}      \left(\int_{B_{2R}}  |x|^\alpha u^{q+1} \right)^{\frac{1}{q+1}}  \left(\int_{B_{2R}}  |x|^{-\frac{\alpha p}{q-p+1}} \right)^{\frac{q-p+1}{p(q+1)}}\\
&=& R^{\frac{n(q-p+1)}{p(q+1)}-\frac{\alpha }{q+1}-1}\left(\int_{B_{2R}}  |\nabla u|^p \right)^{\frac{p-1}{p}}      \left(\int_{B_{2R}}  |x|^\alpha u^{q+1} \right)^{\frac{1}{q+1}}.
\end{eqnarray*}
Therefore, from Step 2, there exists a positive constant $C$ independent of $R$ such that 
\begin{eqnarray*}
|\int_{B_{2R}} |\nabla u|^{p-2} \nabla \zeta_R\cdot \nabla u \ u| & \le & C \ R^{\frac{n(q-p+1)-p(\alpha+q+1)}{p(q+1)}}.
\end{eqnarray*}
Since $n< \frac{p(q+\alpha+1)}{q-p+1}$, we have $\lim_{R\to \infty }|\int_{B_{2R}} |\nabla u|^{p-2} \nabla \zeta_R\cdot \nabla u \ u| =0$. Hence (\ref{identitytest}) implies
\begin{eqnarray}\label{identitypohoz}
\int_{\mathbb{R}^n} |\nabla u|^p = \int_{\mathbb{R}^n}|x|^\alpha u^{q+1}.
\end{eqnarray}

Now, set $\Omega=B_R$ for $R\ge 1$ in (\ref{pohozaevidentity}). Therefore, estimates in Step 2,  i.e.,
$$\int_{0}^{\infty}\int_{|x|=R} |\nabla u|^p< \infty \ \ \text{and} \ \ \int_{0}^{\infty}R^\alpha \int_{|x|=R}u^{q+1}< \infty,$$
imply that R.H.S. of the Pohozaev identity, (\ref{pohozaevidentity}), converges to zero if $R\to \infty$. Hence, 
\begin{eqnarray*} 
\frac{n-p}{p}\int_{\mathbb{R}^n} |\nabla u|^p = \frac{n+\alpha}{q+1}\int_{\mathbb{R}^n} |x|^\alpha u^{q+1}.
\end{eqnarray*}

From this and (\ref{identitypohoz}), we finish the proof of Step 3.

For the second part of the theorem, i.e. the critical case $q=\frac{n(p-1)+p(\alpha+1)}{n-p}, n>p$, the function $ u_{\epsilon} $ defined by (\ref{radial}) satisfies
\begin{equation*}
u_{\epsilon}^{q-1} (|x|) \le C \ \left(\frac{1}{|x|}\right)^{\frac{n(p-2)+p}{p-1}+\frac{p+\alpha}{p-1}} \ \ \ \ \text{and} \ \ \ \ 
|\nabla u_{\epsilon}|^{p-2} (|x|) \ge  C\  \left(\frac{1}{|x|}\right)^{\frac{(n-1)(p-2)}{p-1}},
\end{equation*}
for every $ |x|>R_{0}$, where $R_{0}$ is large enough and the positive constant $C$ does not depend on $|x|$. Since $p+\alpha>0$, for a given $\delta>0$ and $ |x|>R_{0}$, we get 
 \begin{equation*}
u_{\epsilon}^{q-1} (|x|) \le \delta\ \left(\frac{1}{|x|}\right)^{\theta} \ \ \ \ \text{and} \ \ \ \ 
|\nabla u_{\epsilon}|^{p-2} (|x|) \ge C\ \left(\frac{1}{|x|}\right)^{\theta-2},
\end{equation*}
where $\theta:=\frac{n(p-2)+p}{p-1}$. Now, one can apply weighted Hardy's inequality over $\mathbb{R}^{n}\setminus B_{R_{0}}$, i.e., 

$$ \left(\frac{n-\theta}{2}\right)^2  \int_{\mathbb{R}^{n}\setminus B_{R_{0}}} \frac{\phi^2}{|x|^\theta} \le \int_{\mathbb{R}^{n}\setminus B_{R_{0}}} \frac{|\nabla\phi|^2}{|x|^{\theta-2}} ,$$

 to see $u_{\epsilon}$ is stable outside $\overline {B_{R_{0}}}$. Note that $\delta$ can be chosen sufficiently small and for $n>p$ we have $n>\theta$.

\hfill $\Box$

\noindent{\bf Open Problems.} 
\begin{enumerate}
\item The interesting point about Theorem \ref{radialcompact} and \ref{stablenon} is that the higher dimensions given for three different nonlinearities are the same, however the lower dimensions are different.    It would be interesting to see if Theorem \ref{radialcompact} still holds for nonradial solutions and for the same range of parameters.    However, following and adjusting  the same proof, it is straightforward to prove Theorem \ref{radialcompact} for nonradial solutions of (\ref{gener}) in the given dimensions and replacing $\alpha$ by $\alpha^-:=\min\{\alpha,0\}$.    
\item We believe that Theorem \ref{poho1} still holds under the assumption of nonnegative weak solutions and without the extra assumption of finite Morse index.  This seems a challenging problem and even for the semilinear case, that is $p=2$, it is only known in dimension $n=3$ and for bounded solutions by Phan and Souplet in \cite{phs}. A similar result for the H\'{e}non-Lane-Emden system is given in \cite{mostafagh}, Theorem 1.     
\end{enumerate}


\begin{thebibliography}{99}   

 \bibitem{bg} M. F. Bidaut-Veron, H. Giacomini; \emph{A new dynamical approach of Emden-Fowler equations and systems}, Adv. Differential Equations 15 (2010), no. 11-12, 1033-1082.
 
\bibitem{cc} X. Cabr\'e, A. Capella;
 \emph{ On the stability of radial solutions of semi-linear elliptic equations in all of $\mathbb{R}^n$}, C. R. Math. Acad. Sci. Paris 338 (2004) 769-774.

\bibitem{cc1}  X. Cabr\'e, A. Cepella;
 \emph{Regularity of radial minimizers and extremal solutions of semilinear elliptic equations}, J. Funct. Anal. 238 (2006), no. 2, 709-733. 


 
 \bibitem{ccs}  X.  Cabr\'e, A. Capella, M. Sanch\'on;
 \emph{Regularity of radial minimizers of reaction equations involving the $p$-Laplacian}, Calc. Var. Partial Differential Equations 34 (2009), no. 4, 475-494.


\bibitem{cgs} L. Caffarelli, B.  Gidas, J. Spruck;
\emph{Asymptotic symmetry and local behavior of semilinear elliptic equations with critical Sobolev growth}, Comm. Pure Appl. Math. 42 (1989), no. 3, 271-297. 



\bibitem{ces}  D. Castorina, P. Esposito, B.  Sciunzi;
 \emph{Low dimensional instability for semilinear and quasilinear problems in $\mathbb{ R}^n$},  Commun. Pure Appl. Anal.  8  (2009),  no. 6, 1779-1793.



\bibitem{cg} C. Cowan, N. Ghoussoub;
\emph{Estimates on pull-in distances in microelectromechanical systems models and other nonlinear eigenvalue problems}, SIAM J. Math. Anal. 42, no. 5 (2010), 1949-1966.

\bibitem{cf} C. Cowan, M. Fazly;
\emph{On stable entire solutions of semi-linear elliptic equations with weights}, Proc. Amer. Math. Soc. 140 (2012), 2003-2012.




\bibitem{dfsv} L. Damascelli, A. Farina, B. Sciunzi, E. Valdinoci; 
 \emph{Liouville results for $m$-Laplace equations of Lane-Emden-Fowler type},  Ann. Inst. H. Poincar\'e - AN 26 (2009), no. 4, 1099-1119.

%\bibitem {ds} L. Damascelli, B. Sciunzi;
%\emph{Regularity, monotonicity and symmetry of positive solutions of m-Laplace equations}, J. Diff. Eq. 206 (2) (2004), 483-515.


%\bibitem{d2} E. N. Dancer; 
%\emph{Superlinear problems on domains with holes of asymptotic shape and
%exterior problems}, Math. Z., 229 (1998), 475-491.

 \bibitem{df} E. N. Dancer, A. Farina;
\emph{On the classification of solutions of $-\Delta u=e^u$ on $\mathbb{ R}^n$: stability outside a compact set and applications}, Proc. Amer. Math. Soc. 137 (2009), no. 4, 1333-1338.

\bibitem{dg} Y. Du, Z. Guo; 
\emph{ Positive solutions of an elliptic equation with negative exponent: stability and critical power,} J. Diff. Eq. 246 (2009), no. 6, 2387-2414.
 
\bibitem{df1} L. Dupaigne,  A. Farina; 
\emph{Stable solutions of $-\Delta u=f(u)$ in $\mathbb{ R}^n$},  J. Eur. Math. Soc. (JEMS) 12 (2010), no. 4, 855-882. 


\bibitem{e1} P. Esposito; 
\emph{Linear instability of entire solutions for a class of non-autonomous elliptic equations,}  Proc. Roy. Soc. Edinburgh Sect. A  138  (2008),  no. 5, 1005-1018. 



\bibitem{e2} P. Esposito; 
\emph{Compactness of a nonlinear eigenvalue problem with a singular nonlinearity,}  Commun. Contemp. Math.  10  (2008),  no. 1, 17-45. 


 \bibitem{egg1} P. Esposito, N. Ghoussoub, Y. Guo; 
 \emph{ Compactness along the branch of semistable and unstable solutions for an elliptic problem with a singular nonlinearity}. Comm. Pure Appl. Math. 60 (2007), no. 12, 1731-1768.


\bibitem{egg2} P. Esposito, N. Ghoussoub, Y. Guo;
 \emph{Mathematical analysis of partial differential equations modeling electrostatic MEMS}, Courant Lecture Notes in Mathematics, 20. Courant Institute of Mathematical Sciences, New York; American Mathematical Society, Providence, RI, 2010. xiv+318 pp.

\bibitem{f1} A. Farina;
 \emph{On the classification of solutions of the Lane-Emden equation on unbounded domains of $\mathbb{ R}^n$}, J. Math. Pures Appl. (9) 87 (2007), no. 5, 537-561. 



\bibitem{f2} A. Farina;
 \emph{Stable solutions of $-\Delta u=e^u$ on $\mathbb{ R}^n$}, C. R. Math. Acad. Sci. Paris  345 (2007), no. 2, 63-66.


\bibitem{mostafa} M. Fazly;
\emph{Liouville type theorems for stable solutions of 
certain elliptic systems},  Advanced Nonlinear Studies 12 (2012), 1-17.
 
 \bibitem{mostafagh} M. Fazly, N. Ghoussoub;
 \emph{On the H\'{e}non-Lane-Emden conjecture}, To appear Disc. Cont. Dyn. Syst. A.
 
 
 
\bibitem {gg} N. Ghoussoub, Y. J. Guo: 
\emph{On the Partial Differential Equations of Electrostatic MEMS Devices: Stationary Case,} SIAM J. Math. Anal., Vol. 38, No. 5, (2007) p. 1423-1449.

\bibitem{gnn} B. Gidas, W. M. Ni, L. Nirenberg, \emph{Symmetries and related properties via the maximum
principle}, Comm. Math. Phys. 68 (1979), 209-243.

%\bibitem{gy} N. Ghoussoub, C. Yuan; \emph{Multiple Solutions for Quasi-linear PDEs Involving the Critical Sobolev and Hardy Exponents}, Trans. Amer. Math. Soc. 352 (2000), 5703-5743.


\bibitem{gs1} B. Gidas,  J. Spruck;
\emph{Global and local behavior of positive solutions of nonlinear
elliptic equations}, Commun. Pure Appl. Math. 34 (1981) 525-598.

\bibitem {gs2} B. Gidas, J. Spruck;
\emph{A priori bounds for positive solutions of nonlinear elliptic equations,} Comm. Partial Differential Equations 6 (1981), no. 8, 883-901.
 





\bibitem{gt} D. Gilbarg, N. S. Trudinger;
\emph{Elliptic partial differential equations of second order}, Springer-Verlag, Berlin, 2001. xiv+517 pp.


\bibitem{gpw} Y. Guo, Z. Pan, M. J. Ward;
\emph{Touchdown and pull-in voltage behavior of a MEMS device
with varying dielectric properties}, SIAM J. Appl. Math., 66 (2005), 309-338.

\bibitem{h} M. H\'{e}non, 
\emph{Numerical experiments on the stability of spherical stellar systems,} Astron. Astrophys. 24 (1973) 229-238.

\bibitem{n} W. M. Ni;  
\emph{A nonlinear Dirichlet problem on the unit ball and its applications,}
Indiana Univ. Math. J. 31 (1982), no. 6, 801-807.
 
%\bibitem{ps} A. L. Peletier, J. Serrin; 
%\emph{Gradient bounds and Liouville theorems for quasilinear
%elliptic equations},  Ann. Scuola Norm. Sup. Pisa C1. Sei (4), 5 (1978), 65-104.

\bibitem{pb}  J. Pelesko, D. Bernstein;   \emph{Modeling MEMS and NEMS}, Chapman $\&$ Hall/CRC, Boca Raton, FL, 2003. xxiv+357 pp. 


\bibitem{phs} Q. H. Phan, Ph. Souplet; 
\emph{Liouville-type theorems and bounds of solutions of Hardy-H\'{e}non equations}, J. Differential Equations 252 (2012), no. 3, 2544-2562. 


%\bibitem {pqs} P. Pol\'{a}\v{c}ik, P. Quittner, Ph. Souplet; \emph{Singularity and decay estimates in superlinear problems via Liouville-type theorems,} Part I: Elliptic equations and systems, Duke Math. J. 139 (2007) 555-579.

\bibitem {s1} J. Serrin; 
\emph{Local behavior of solutions of quasi-linear elliptic equations}, Acta Math. 111 (1964)
247-302.

%\bibitem{s2} J. Serrin;
%\emph{Entire solutions of nonlinear Poisson equations}, Proc. London Math. Soc. (3), 24 (1972),
%348-366.

\bibitem{sz} J. Serrin, H. Zou;
\emph{Cauchy-Liouville and universal boundedness theorems for quasilinear elliptic equations and inequalities }, Acta Mathematica, 189 (2002) 79-142.

\bibitem{ssw} D. Smets, J. Su, M. Willem, \emph{Non radial ground states for the H\'{e}non equation}, Commun.
Contemp. Math. 4 (2002), no. 3, 467-480.

\bibitem{sv}  S. Villegas;
\emph{Asymptotic behavior of stable radial solutions of semilinear elliptic equations in $\mathbb{R}^n$}, J. Math. Pures Appl. (9) 88 (2007), no. 3, 241-250.

% \bibitem {w} X. Wang,  \emph{On the Cauchy Problem for Reaction-Diffusion Equations,} Tran. Amer. Math. Soc., Vol. 337, No. 2 (1993), 549-590.

\bibitem{wy} C. Wang, D. Ye,  \emph{Some Liouville theorems for H\'{e}non type elliptic equations},  J. Funct. Anal. 262 (2012), no. 4, 1705-1727.
  
\end{thebibliography}
\end{document}